\input amstex
\input amsppt.sty   
\hsize 26pc
\vsize 46pc
\magnification=\magstep1
\def\nmb#1#2{#2}         
\def\cit#1#2{\ifx#1!\cite{#2}\else#2\fi} 
\def\totoc{}             
\def\idx{}               
\def\ign#1{}             

\redefine\o{\circ}

\define\al{\alpha}

\define\de{\delta}
\define\ep{\varepsilon}

\define\ka{\kappa}
\define\la{\lambda}
\define\rh{\rho}
\define\si{\sigma}
\define\ta{\tau}
\define\ph{\varphi}

\define\De{\Delta}

\define\La{\Lambda}

\define\row#1#2#3{#1_{#2},\ldots,#1_{#3}}
\define\rowup#1#2#3{#1^{#2},\ldots,#1^{#3}}
\define\x{\times}
\define\Lam#1#2#3{\La^{#1}(#2;#3)}
\redefine\S{{\Cal S}}
\define\A{{\Cal A}}
\define\s#1#2{\operatorname{sign}(#1,\bold #2)}
\define\Der{\operatorname{Der}}
\define\End{\operatorname{End}}
\def\today{\ifcase\month\or
 January\or February\or March\or April\or May\or June\or
 July\or August\or September\or October\or November\or
December\fi \space\number\day, \number\year}
\topmatter
\title The multigraded Nijenhuis-Richardson Algebra,\\ 
its universal Property and Applications   \endtitle
\author P.A.B. Lecomte \\ 
        P.W. Michor \\ 
        H. Schicketanz \endauthor
\affil Institut f\"ur Mathematik der Universit\"at Wien, Austria \\
Institut des Math\'ematique, Universit\'e de Li\`ege, Belgium
\endaffil
\address{P\. A\. B\. Lecomte: 
Universit\'e de Li\`ege, Institut de Math\'ematique, Avenue
des Tilleuls, 15; B-4000 Li\`ege, Belgium}\endaddress
\address{P\. W\. Michor, H\. Schicketanz: 
Institut f\"ur Mathematik, Universit\"at Wien,
Strudlhofgasse 4, A-1090 Wien, Austria}\endaddress
\date{February 22, 1991
}\enddate
\abstract{
We define two $(n+1)$ graded Lie brackets on spaces of multilinear 
mappings. The first one is able to recognize $n$-graded associative 
algebras and their modules and gives immediately the correct 
differential for Hochschild cohomology. The second one recognizes 
$n$-graded Lie algebra structures and their modules and gives rise to 
the notion of Chevalley cohomology.}\endabstract
\endtopmatter

\leftheadtext{\smc Lecomte, Michor, Schicketanz}
\rightheadtext{\smc Multigraded Nijenhuis Richardson algebra}
\document

\heading \nmb0{1}. Introduction \endheading

In this paper we will generalize the construction of
Nijenhuis and Richardson which associates to a given vector space
$V$ a graded Lie algebra $Alt(V)$ of multilinear alternating 
mappings $V\x\dots\x V \to V$ to study Lie algebra
structures on $V$ and their deformations, see \cit!{9}. Their
construction suggests a "principle" which we present here as the
starting point for our investigations. The principle is as
follows:

Suppose that $\S$ is a type of structures on $V$, defining for
example associative algebras, Lie algebras, modules (over a
given Algebra $\A$) or Lie bialgebras on $V$. Then there exists
a $\Bbb Z$ - graded Lie algebra 
$( \Cal E = \bigoplus_{k \in  \Bbb Z}\Cal E^k, [\quad,\quad] )$ 
such that $P \in \S$ if and only if $P \in \Cal E^1$
such that $[P,P] = 0$.

In the case where $\S$ is the set of Lie algebra structures on
$V$ the space $\Cal E$ can be identified with $Alt(V)$. Moreover if
$V$ is equipped with such a $P$, the Chevalley-Eilenberg
coboundary operator $\partial_P$ of the adjoint representation of
$(V,P)$ is just the adjoint action of $P$ on
$Alt(V)$  up to a sign.  
Another application may be found in \cit!{6}. There
one uses the graded cohomology of the subalgebras of $Alt(V)$ to
classify and to construct formal deformations of $(V,P)$.

The purpose of this paper is to establish the principle in each of
the cited cases. We will do this in more generality which makes
the construction even more powerful. Namely, we assume that $V$ is itself
graded over $\Bbb Z^n$, $(n = 0,1,2,...)$ and we will define for
each $\S$ a graded Lie algebra $\Cal E$ which is now graded over
$\Bbb Z^{n+1}$ and satisfies the principle. If we don't stress
the special choice of $n$ we will speak of {\it multigraded}
algebras. Having defined the multigraded Lie algebra $\Cal E$, deformation 
theory and cohomology of $\S$ may be treated at the same time using only
the space $\Cal E$ and its properties.

Given a multigraded vector space $V$ we will construct first
$M(V)$, a multigraded Lie algebra which is adapted to study the
associative structures on $V$. Using then the multigraded 
alternator $\al$ we define $A(V)$ to be the image of $M(V)$ by
$\al$ equipped with the unique bracket making $\al$ a
homomorphism of multigraded Lie algebras. Moreover $A(V)$
satisfies a universal property and describes multigraded Lie
algebra structures on $V$. We call $A(V)$ the {\it multigraded
Nijenhuis-Richardson algebra} of $V$ since it coincides with
$Alt(V)$ for $n = 0$. Once having established
this multigraded version, the result for module
structures follows quite easily.

In this way we rediscover Hochschild and Chevalley-Eilenberg
Cohomology for $n \le 1$, where the differential is given by the
adjoint action of $P$ on $\Cal E$. Their generalizations for 
$n>1$ are now obvious and yield a canonical description for
multigraded cohomology in both cases. 

Moreover one can study now the theory of formal deformations of 
multigraded algebras $L$ and their modules. Roughly speaking we 
describe a mapping from the cohomology of the adjoint representation 
of $A(L)$ into the set of formal deformations of all possible 
structures on $L$ which may be used to construct and classify these 
deformations. Such a point of view has also been emphasized by 
\cit!{11}, \cit!{5}, and \cit!{4}.

\heading \nmb0{2}.Multigraded associative algebra structures \endheading

\subheading{\nmb.{2.1}. Conventions and definitions.} By a {\it
multidegree} we mean an element $x = (\rowup x1n) \in \Bbb Z^n$
for some $n$. We call it also {\it $n$-degree} if we want to
stress the special choice of $n$. We shall need also the inner
product of multidegrees $\langle \quad,\quad\rangle : \Bbb Z^n \x \Bbb Z^n \to
\Bbb Z$, given by $\langle x,y\rangle  = \sum_{i=1}^n x^iy^i$.

An {\it $n$-graded vector space} is just a direct sum $V =
\bigoplus_{x \in \Bbb Z^n}V^x$, where the elements of $V^x$ are
said to be homogeneous of multidegree $x$. To avoid technical
problems we assume that vector spaces are defined
over a field $\Bbb K$ of characteristic 0. In the following $X$,
$Y$, etc will always denote homogeneous elements of some
multigraded vector space of multidegrees $x$, $y$, etc.

By an $n$-graded algebra $\A  = \bigoplus_{x \in
\Bbb Z^n}\A ^x$ we mean an $n$-graded vector space which is also
a $\Bbb K$ algebra such that $\A^x\cdot \A^y \subseteq \A^{x+y}$.
\roster \item The multigraded algebra $(\A,\cdot)$ is said to be {\it
multigraded commutative} if for homogeneous elements $X$, $Y \in
\A $ of multidegree $x$, $y$, respectively,we have $X\cdot Y =
(-1)^{\langle x,y\rangle }Y\cdot X$.
\item If $X\cdot Y = -(-1)^{\langle x,y\rangle }Y\cdot X$ holds it is 
     called {\it multigraded       anticommutative}.
\item An {\it $n$-graded Lie algebra} is a multigraded
     anticommutative algebra $(\Cal E,[\quad,\quad])$, such that the {\it
     multigraded Jacobi identity} holds: 
     $$[X,[Y,Z]] = [[X,Y],Z] + (-1)^{\langle x,y\rangle }[Y,[X,Z]]$$
\endroster
Obviously the space 
$\End(V) = \bigoplus_{\de \in \Bbb Z ^n}\End^{\de}(V)$
of all endomorphisms of a multigraded vector space
$V$ is a multigraded algebra under composition, where
$\End^\de(V)$ is the space of linear endomorphisms $D$ of $V$ of
multidegree $\de$, i.e. $D(V^x) \subseteq V^{x+\de}$. Clearly
$\End(V)$ is a multigraded Lie algebra under the multigraded
commutator
$$[D_1,D_2] := D_1\circ D_2 + (-1)^{\langle \de_1,\de_2\rangle } D_2\circ D_1.\tag4$$

If $\A $ is an $n$-graded algebra, an endomorphism 
$D:\A \to \A$ of multidegree $\de$ is called a {\it multigraded
derivation}, if for $X$, $Y \in \A $ we have 
$$D(X\cdot Y) = D(X) \cdot Y + (-1)^{\langle \delta,x\rangle }X\cdot D(Y).\tag5$$ 
Let us write $\Der^\de (\A )$ for the space of all multigraded
derivations of degree $\de$ of the algebra $\A $, and we put 
$$\Der(\A ) = \bigoplus_{\de \in \Bbb Z^n}\Der^\de (\A ).\tag5$$

The following lemma is standard:
\proclaim{Lemma} If $\A $ is an $n$-graded algebra, then the space
$\Der(\A )$ of multigraded derivations is an $n$-graded Lie 
subalgebra under the $n$-graded commutator.
\endproclaim

It is clear from the definitions that non-graded algebras and
$\Bbb Z$-graded algebras are multigraded of multidegree 0 and 
1, respectively.

\subheading{\nmb.{2.2} Associative algebra structures} Let us
recall first the construction in the case of non-graded vector
spaces which was given in \cit!{3}, \cit!{1}. 
There a $1$-graded Lie algebra $(M(V),[\quad,\quad]^\De)$ is 
described for each vector space $V$ with 
the property that $(V,\mu)$ is an associative algebra if and only if 
$\mu \in M^1(V)$ and $[\mu,\mu]^{\De}=0$. This algebra is as follows.

Denote by $M^k(V)$ the space of all  $k+1$-linear mappings
$K : V\x \dots \x V \to V$ 
and set
$$M(V):=\bigoplus_{k\in \Bbb Z}M^k(V).$$
For $K_i\in M^{k_i}(V)$ and $X_j\in V$ we define 
$j(K_1)K_2 \in M^{k_1+k_2}(M)$ by
$$\align
&(j(K_1)K_2)(\row X0{k_1+k_2}) := \\
&\quad =\sum_{i=0}^{k_2}(-1)^{k_1i}
    K_2(X_0,\dots,K_1(\row Xi{i+k_1}),\dots,X_{k_1+k_2}).
\endalign $$
The graded Lie bracket of $M(V)$ is then given by
$$[K_1,K_2]^\De=j(K_1)K_2-(-1)^{k_1k_2}j(K_2)K_1.$$

\proclaim{Proposition} {\rm (\cit!{3}, \cit!{1})}
\roster
\item $(M(V),[\quad,\quad]^\De)$ is a 1-graded Lie algebra.
\item If $\mu\in M^1(V)$, so $\mu:V\x V\to V$ is bilinear, then 
    $(V,\mu)$ is an associative algebra if and only if 
    $[\mu,\mu]^\De=0$. \qed
\endroster
\endproclaim

Note that $M^0(V)=\End(V)$ is a Lie subalgebra of $M(V)$, and its 
bracket is the negative of the usual commutator.

The explicit formulas above follow directly from investigating the 
1-graded Lie algebra of (1-graded) derivations of certain graded 
algebras, see \cit!{11}. We explain that in the simple case of a 
finite dimensional $V$.
Then $M(V)$ is canonically isomorphic 
to the $1$-graded Lie algebra $\Der(\bigotimes V^*)$ of derivations 
of the tensor algebra of $V^*$, a derivation $D$ of degree $k$ being 
completely determined by its restriction $V*\to \bigotimes^{k+1}V^*$ 
and hence by a unique $K\in M^k(V)$.

\subheading{\nmb.{2.3} Multigraded associative algebras} 
We will give now the multigraded generalization. Of course on can 
proceed as before by identifying $M(V)$ as the algebra of derivations
of some suitable multigraded algebra. But we will generalize 
\nmb!{2.2} directly. So let $V=\bigoplus_{x\in\Bbb Z^n}V^x$ be 
an $n$-graded vector space. We define 
$$M(V):=\bigoplus_{(k,\ka)\in\Bbb Z\x\Bbb Z^n}M^{(k,\ka)}(V),$$
where $M^{(k,\ka)}(V)$ is the space of all $k+1$-linear mappings
$K:V\x\dots\x V\to V$ such that 
$K(V^{x_0}\x\dots\x V^{x_k})\subseteq V^{x_0+\dots+x_k+\ka}$.
We call $k$ the \idx{\it form degree} and $\ka$ the \idx{\it weight 
degree} of $K$. In \nmb!{2.2} the mapping 
$K$ had degree $k$ and $X_i$ had degree 
$-1$ in $M(V)$, hence the sign $(-1)^{ki}$.
We define for $K_i\in M^{(k_i,\ka_i)}(V)$ and $X_j\in V^{x_j}$
$$\multline(j(K_1)K_2)(\row X0{k_1+k_2}) := \\
=\sum_{i=0}^{k_2}
(-1)^{k_1i+\langle \ka_1,\ka_2+x_0+\dots+x_{i-1}\rangle}\cdot\\
\shoveright{\cdot K_2(X_0,\dots,K_1(\row Xi{i+k_1}),\dots,X_{k_1+k_2})}\\
\shoveleft{[K_1,K_2]^\De=
j(K_1)K_2-(-1)^{k_1k_2+\langle \ka_1,\ka_2\rangle}j(K_2)K_1.\hfill}
\endmultline$$

\proclaim{Theorem} Let $V$ be an $n$-graded vector
space. Then we have:
\roster
\item $(M(V),[\quad,\quad]^\De)$ is an $(n+1)$-graded Lie algebra.
\item If $\mu\in M^{(1,0,\dots,0)}(V)$, so $\mu:V\x V\to V$ is bilinear 
of weight $0\in\Bbb Z^n$, then $\mu$ is an associative $n$-graded 
multiplication if and only if $[\mu,\mu]^\De=0$.
\endroster
\endproclaim

\demo{Proof}
The bracket is $(n+1)$-graded anticommutative. The $(n+1)$-graded 
Jacobi identity follows from the formula
$$j([K_1,K_2]^\De)=[j(K_1),j(K_2)],$$
the multigraded commutator in $\End(M(V))$. This is a long but elementary 
calculation. The second assertion follows by writing out the definitions.
\qed\enddemo

\heading\totoc\nmb0{3}. Multigraded Lie Algebra Structures
\endheading

\subheading{\nmb.{3.1}. Multigraded signs of permutations} 
Let $\bold x = (\row x1k) \in (\Bbb Z^n)^k$ be a multi index of $n$-degrees
$x_i = (x_i^1,\dots,x_i^n) \in \Bbb Z^n$ and let $\si \in \S_k$ be a
permutation of $k$ symbols. Then we define the {\it multigraded
sign} $\s \si x$ as follows: 
For a transposition $\si =(i,i+1)$ we
put $\s \si x = -(-1)^{\langle x_i,x_{i+1}\rangle }$; 
it can be checked by combinatorics that this gives a well defined
mapping $\s {\quad}x:\Cal S_k \to  \{-1,+1\}$. In fact one may define
directly
$$\s \si x = \operatorname{sign}(\si)
\operatorname{sign}(\si_{|x_1^1|,\dots,|x_k^1|})\cdots
\operatorname{sign}(\si_{|x_1^n|,\ldots,|x_k^n|}),$$
where $\si_{|x_1^j|,\ldots,|x_k^j|}$ is that permutation of 
$|x_1^j|+\dots+|x_k^j|$ symbols which moves the $i$-th block of
length $|x_j^i|$ to the position $\si i$, and where 
$\operatorname{sign}(\si)$ 
denotes the ordinary sign of a permutation in $\S_k$.
Let us write $\si x = (\row x{\si1}{\si k})$, then we have the following 

\proclaim{Lemma} $\s {\si\o\ta}x = \s \si x .\s \ta{\si \bold x}$.\qed
\endproclaim

\subheading{\nmb.{3.2} Multigraded Nijenhuis-Richardson algebra}
We define the {\it multigraded alternator} $\al : M(V)\to M(V)$ by
$$(\al K)(X_0,\dots,X_k) = \dsize\frac1{(k+1)!} \sum_{\si \in \S_{k+1}}
\s \si x K(X_{\si 0},\dots,X_{\si k})\tag1$$
for $K \in M^{(k,*)}(V)$ and $X_i \in V^{x_i}$. 
If the ground field is not of characteristic 0 one could omit the 
combinatorial factor, but one should redo the whole developpment
starting from the point of view of derivations again, see the remark 
at the end of \nmb!{2.2}. However, the combinatorial factors used here are 
quite essential, judging from our experience in differential geometry.
By lemma \nmb!{3.1} we have
$\al^2 = \al$ so $\al$ is a projection defined on $M(V)$, homogeneous of 
multidegree 0, and we set 
$$\align
A(V)&=\bigoplus_{(k,\ka)\in\Bbb Z\x \Bbb Z^n}A^{(k,\ka)}(V)\\
:&=\bigoplus_{(k,\ka)\in \Bbb Z\x \Bbb Z^n}\al (M^{(k,\ka)}(V)).
\endalign$$
A long but straightforward computation shows that for 
$K_i \in M^{(k_i,\ka_i)}(V)$ 
$$\al (j(\al K_1)\al K_2) = \al (j(K_1)K_2),$$
so the following operator and bracket is well defined:
$$\align
i(K_1)K_2 :&= \dsize\frac{(k_1+k_2+1)!}{(k_1+1)!(k_2+1)!}
     \al(j(K_1)K_2)\\
[K_1,K_2]^{\wedge} &=\dsize\frac{(k_1+k_2+1)!}{(k_1+1)!(k_2+1)!}
                                               \al ([K_1,K_2]^{\De}) \\
&=i(K_1)K_2 -(-1)^{\langle (k_1\kappa_1),(k_2,\kappa_2)\rangle }i(K_2)K_1
\endalign$$
The combinatorial factor will become clear in \nmb!{3.4} .

\proclaim{\nmb.{3.3}. Theorem}
1. If $K_i$ are as above then
$$\multline
(i(K_1)K_2)(\row X0{k_1+k_2}) = \\
=\dsize\frac1{(k_1+1)!k_2!}\sum_{\si \in \S_{k_1+k_2+1}}\s \si x 
     (-1)^{\langle \ka_1,\ka_2\rangle }\cdot\\
\cdot K_2((K_1(X_{\si 0},\dots,X_{\si k_1}),\dots,X_{\si(k_1+k_2)}).
\endmultline$$

2. $(A(V),[\quad,\quad]^{\wedge})$ is an $(n+1)$-graded Lie algebra.

3. If $\mu\in A^{(1,0,\dots,0)}(V)$, so $\mu:V\x V\to V$ is bilinear 
$n$-graded anticommutative mapping of weight $0 \in \Bbb Z^n$ then 
$[\mu ,\mu ]^{\wedge}=0$ if and only if $(V,\mu)$ is a $n$-graded Lie algebra.
\endproclaim

\demo{Proof}
1. This follows by a straight forward computation.

2. $[\quad,\quad]^{\wedge}$ is clearly multigraded anticommutative 
and the	multigraded Jacobi identity follows directly from the one of 
$[\quad,\quad]^\De$.

3. Let $\mu \in A^{(1,0,\dots,0)}(V)$, then
$$\align
0 &= [\mu,\mu]^{\wedge} (X_0,X_1,X_2) \\
&= \dsize \frac{3!}{2!3!}\frac1{3!}\sum_{\si \in \Cal S_3}\s \si x\cdot
      [\mu,\mu]^{\De}(X_{\si0},X_{\si1},X_{\si2})\\
&=\sum_{\si \in \Cal S_3}\s \si x\cdot \mu(\mu(X_{\si0},
      \mu( X_{\si1},X_{\si2})) \\
\endalign$$
which is equivalent to the multigraded Jacobi identity of $(V,\mu).$
\qed\enddemo

We call $(A(V),[\quad,\quad]^{\wedge})$ the {\it multigraded 
Nijenhuis-Richardson algebra}, since $A(V)$ coincides for $n=0$ 
with $Alt(V)$ of \cit!{9}.

\subheading{\nmb.{3.4}. Universality of the algebra 
$(A(V),[\quad,\quad]^{\wedge})$} Let $V$ be a multigraded vector space 
and denote by $\Cal E (V)$ the category of multigraded Lie algebras 
$(E,[\quad,\quad])$ such that
$$\align
\quad \quad E^{(k,*)} &= 0 \quad k<-1 \\
	    E^{(-1,*)}&= V .
\endalign$$
If $E,F \in \Cal E (V)$, then a {\it morphism} $\ph : E \to F$ is a 
homomorphism of multigraded Lie algebras satisfying 
$\ph | E^{(-1,*)}= \operatorname{id}_V$. 
For example $M(V)$ and $A(V)$ are elements of $\Cal E(V)$.

\proclaim{Theorem}
$A(V)$ is a final object in $\Cal E (V)$, so for each $E \in \Cal E(V)$ 
there exists a unique morphism $\ep : E \to A(V)$. 
It follows that $A(V)$ is unique up to isomorphism.
\endproclaim

\demo{Proof}
Suppose that $Z \in E^{(k,z)}$ then we define
$$\ep (Z) (\row X0k) = (-1)^{\langle z,x_0+\dots +x_k\rangle }
     [X_0,[X_1,\dots ,[X_k,Z]\dots],$$
an element of $E^{(-1,*)}=V$ for 
$X_i \in V^{x_i}$. Because of the multigraded Jacobi identity 
$\ep (Z)$ is well defined as an element of $A^{(k,z)}$. So we are 
left to	show that
$$\ep ([Z_1,Z_2]) = [\ep (Z_1),\ep (Z_2)]^{\wedge}\tag*$$
We will do this by induction on $k=k_1+k_2$. For $k<-1$ this is 
trivially
true. Now let $k=-1$, so we may assume that $Z_1 \in V^{z_1}$. Then
$$\align
\ep ([Z_1,Z_2] &= [Z_1,Z_2] = (-1)^{\langle z_1,z_2\rangle }\ep (Z_2)(Z_1)\\
&=i(Z_1)\ep (Z_2) = [Z_1,\ep (Z_2)]^{\wedge} 
     = [\ep (Z_1),\ep (Z_2)]^{\wedge} 
\endalign$$
by {Theorem 3.2} and since $\ep | V = \operatorname{id}_V$.
Suppose that \thetag* is true for $k_1+k_2 < k$. Then for $k_1+k_2 = 
k$ we have
$$\align
&i(X)\ep ([Z_1,Z_2])=[X,\ep ([Z_1,Z_2])]^{\wedge} 
=\ep ([X,[Z_1,Z_2]]) \\
&\quad =\ep \bigl([[X,Z_1],Z_2] + 
(-1)^{\langle (-1,x),(k_1,z_1)\rangle }[Z_1,[X,Z_2]]\bigr)  \\
&\quad =[\ep ([X,Z_1]),\ep (Z_2)]^{\wedge} 
     + (-1)^{\langle (-1,x),(k_1,z_1)\rangle }[\ep (Z_1),\ep([X,Z_2])]^{\wedge} \\
&\quad =[i(X)\ep (Z_1),\ep (Z_2)]^{\wedge} 
     + (-1)^{\langle (-1,x),(k_1,z_1)\rangle }[\ep (Z_1),i(X)\ep (Z_2)]^{\wedge} \\
&\quad =i(X)[\ep (Z_1),Z_2]^{\wedge}
\endalign$$
by induction hypothesis and the fact that $i(X)= [X,\quad]^{\wedge}$ is a 
derivation of degree $(-1,x)$ of $A(V)$. This proves the induction. Remark 
that for $E=M(V)$ the morphism $\ep$ is given by
$$\ep | {M^{k,*}(V)} = (k+1)!\;\al.\qed$$
\enddemo

\heading \nmb0{4}. Multigraded Modules and Cohomology \endheading

\subheading{\nmb.{4.1}. Multigraded bimodules}
Let $V$ and $W$ be multigraded vector spaces and $\mu:V\x V \to V$ a 
multigraded algebra structure. A {\it multigraded bimodule} 
$\Cal M=(W,\la,\rh)$ over $\A =(V,\mu)$
is given by $\la ,\rh : V\to \End(W)$ of weight $0$ such that
$$\align
[\mu,\mu]^{\De} &= 0 \quad \text{ so } \A \text{ is associative } \tag1 \\
\la (\mu (X_1,X_2)) &= \la (X_1)\circ \la (X_2)  \tag2 \\
\rho (\mu (X_1,X_2)) &= (-1)^{\langle x_1,x_2\rangle}
     \rho (X_2)\circ \rho (X_1)  \tag3 \\
\la (X_1)\circ \rho (X_2) &= (-1)^{\langle x_1,x_2\rangle }
     \rho (X_2)\circ \la (X_1) \tag4 
\endalign$$
where $X_i \in V^{x_i}$ and $\circ$ denotes the composition in $\End (W)$.

\proclaim{\nmb.{4.2}. Theorem}
Let $E$ be the multigraded vector space defined by
$$E^{(k,*)}=\left\{\alignedat2 &V &\qquad &\text{if } k=0 \\ 
                    &W &&\text{if } k=1 \\ 
                    &0 && \text{otherwise.} 
               \endalignedat
\right. $$
Then $P\in M^{(1,0,\dots,0)}(E)$ defines a bimodule structure on $W$ 
if and only if $[P,P]^{\De}=0$.
\endproclaim

\demo{Proof}
We define
$$\align
\mu (X_1,X_2)&:=P(X_1,X_2) \\
\la (X)Y&:=P(X,Y) \\
\rho(X)Y&:=(-1)^{\langle x,y\rangle }P(Y,X)
\endalign$$
where we suppose the $X_i$'s $\in V$ and $Y \in W$ to be embedded in $E$.
Then if $Z_i \in E$ be arbitrary we get
$$\align
[P,P]^\De(Z_0,Z_1,Z_2) & =2(j(P)P)(Z_0,Z_1,Z_2)	  \\
&= 2P((Z_0,Z_1),Z_2) - 2P(Z_0,(Z_1,Z_2)).
\endalign$$
Now specify $Z_i\in V$ resp. $W$ to get eight independent equations. Four
of them vanish identically because of their degree of homogeneity, the 
others recover the defining equations for the multigraded bimodules.
\qed\enddemo

\proclaim{\nmb.{4.3} Corollary}
In the above situation we have the following decomposition of $M(E)$ :
$$M^{(k,q,*)}(E) = \left\{
     \alignedat2 &0 &\quad &\text{for } q>1 \\
          &L^{(k+1,*)}(V,W) &&\text{for } q=1 \\
          &M^{(k,*)}(V)\oplus\bigoplus^{k+1}(L^{(k,*)}(V,\End(W)) &&\text{for } q=0
     \endalignedat
\right.$$
where $L^{(k,*)}(V,W)$ denotes the space of $k$-linear mappings 
$V\x\dots\x V\to W$.
If $P$ is as above, then
$P=\mu +\la +\rho$ corresponds exactly to this decomposition.
\qed\endproclaim

\subheading{\nmb.{4.4}. Hochschild cohomology and multiplicative 
structures}
Let $V$,$W$ and $P$ be as in {Theorem \nmb!{4.2}} and let $\nu : W\x W \to W$ be
a multigraded algebra structure, so $\nu \in M^{(1,-1,0,\dots,0)}(E)$.
Then for $C_i\in L^{(k_i,c_i)}(V,W)$ we define
$$C_1\bullet C_2:=[C_1,[C_2,\nu]^{\De}]^{\De}.$$
Since $[C_1,D_2]^{\De}=0$ it follows that $(L(V,W),\bullet)$ is multigraded
commutative. It is the usual extension of the product $\nu$ from $W$ 
to the level of cochains, where the necessary combinatorics is hidden 
in the brackets.

\proclaim{Theorem}
1. The mapping $[P,\quad]^{\De} : M(E)\to M(E)$ is a differential. We 
denote its restriction to $L(V,W)$ by $\de_P$. This generalizes
the Hochschild coboundary operator to the multigraded case: If
$C \in L^{(k,c)}(V,W)$ then we have for $X_i \in V^{x_i}$
$$\multline 
(\de_P C)(\row X0k) = \la (X_0)C(\row X1k) \\
-\sum_{i=0}^{k-1}(-1)^i C(X_0,\dots,\mu (X_i,X_{i+1}),\dots,X_k) \\
+(-1)^{k+1+\langle x_0+\cdots+x_{k-1}+c,x_k\rangle}\rho (X_k)C(\row X0{k-1})
\endmultline$$
The corresponding cohomology will be denoted by $H(\A,\Cal M)$, where
$\Cal A$ is the multigraded associative algebra $(V,\mu)$, and where 
$\Cal M$ is the multigraded $\Cal A$-bimodule $(W,\la,\rh)$

2. If $[P,\nu]^{\De}=0$ then $\de_P$ is a derivation of $L(V,W)$ of 
multidegree $(1,0,...0)$. In this case the product $\bullet$ carries over
to a multigraded (cup) product on $H(\A,\Cal M)$.
\endproclaim

\demo{Proof}
The fact that $\de_P$ is a differential follows directly from the 
multigraded Jacobi identity since the degree of $\de_P$ is 
$(1,0,\dots,0)$. The formula is easily checked by writing 
out the definitions. Applying the multigraded Jacobi identity 
once again one gets immediately that $\de_P$ is a derivation 
if and only if $[P,\nu]^{\De}=0$. 

By writing out the definitions one shows that $[P,\nu]^{\De}=0$ is
equivalent to the following equations:
$$\align
\la (X)\nu (Y_1,Y_2)&=\nu (\la (X)Y_1,Y_2)) \\
\rho (X)\nu (Y_1,Y_2)&=(-1)^{\langle x,y_1\rangle }\nu (Y_1,\rho (X)Y_2) \\
\nu (\rho (X)Y_1,Y_2)&=(-1)^{\langle x,y_1\rangle }\nu (Y_1,\la (X)Y_2) \\
\endalign$$
in particular we have $(\la -\rho) : V \to \Der(W,\nu)$.
\qed\enddemo

\subheading{\nmb.{4.5} Multigraded Lie modules and Chevalley cohomology}
We obtain a corresponding result for Lie modules by applying the multigraded
alternator $\al$ to $M(E)$, just as we did in section \nmb!{3} to obtain the 
Nijenhuis-Richardson bracket.

\proclaim{Theorem}
Let $P\in A^{(1,0,\dots,0)}(E)$ then $[P,P]^{\wedge}=0$ if and
only if
$$[\mu,\mu]^{\wedge}=0\tag a$$  
so $(V,\mu)=\frak g$ is a multigraded Lie algebra, and 
$$\pi (\mu (X_1,X_2))Y = [\pi (X_1),\pi (X_2)]Y \tag b$$
where $\mu (X_1,X_2)=P(X_1,X_2)\in V$ and $\pi (X)Y=P(X,Y)\in W$
for $X$, $X_i\in V$ and $Y\in W$, and where $[\quad,\quad]$
denotes the multigraded commutator in $\End (W)$. So $[P,P]^{\wedge}=0$
is by definition equivalent to the fact that $\Cal M :=(W,\pi)$ is a 
multigraded Lie-$\frak g$ module.

If $P$ is as above the mapping 
$\partial_P := [P,\quad]^{\wedge}: A(E) \to A(E)$ is a differential 
and its restriction to 
$$\bigoplus_{k\in \Bbb Z}
\La^{(k,*)}(\frak g,\Cal M):= \bigoplus_{k\in \Bbb Z}A^{(k,1,*)}(E)$$
generalizes the  Chevalley-Eilenberg coboundary operator to the 
multigraded case:
$$\align
(\partial_PC)(\row X0k) 
     &= \sum_{i=0}^k(-1)^{\al_i(\bold x)+\langle x_i,c\rangle }
     \pi (X_i)C(X_0,\dots,\widehat{X_i},\dots,X_k) \\
&+\sum_{i<j}(-1)^{\al_{ij}(\bold x)}
     C(\mu (X_i,X_j),\dots,\widehat{X_i},\dots,\widehat{X_j},\dotsc)
\endalign$$
where 
$$\left\{
\aligned \al_i(\bold x)&=\langle x_i,x_1+\dots+x_{i-1}\rangle +i \\
\al_{ij}(\bold x)&=\al_i(\bold x)+\al_i(\bold x)+\langle x_i,x_j\rangle 
\endaligned\right.$$
We denote the corresponding cohomology space by $H(\frak g,\Cal M)$.

If $\nu : W\x W\to W$ is multigraded symmetric (so 
$\nu \in A^{(1,-1,*)}(E)$) and $[P,\nu]^{\wedge}=0$ then $\partial_P$ acts
as derivation of multidegree $(1,0,\dots,0)$ on the multigraded commutative
algebra $(\La (\frak g,\Cal M),\bullet)$, where
$$C_1\bullet C_2:=[C_1,[C_2,\nu]^{\wedge}]^{\wedge}
     \quad C_i\in \La^{(k_i,c_i)}(\frak g,\Cal M).$$
In this situation the product $\bullet$ carries over 
to a multigraded symmetric
(cup) product on $H(\frak g,\Cal M)$.
\endproclaim

\demo{Proof}
Apply the multigraded alternator $\al$ to the results of \nmb!{4.1}, 
\nmb!{4.2}, \nmb!{4.3}, and \nmb!{4.4}.
\qed
\enddemo

The formulas we obtained here are not that surprising since they are
standard in the non-graded case. The new feature of our approach lies
in the fact that we can formulate deformation equations and cohomology
at once inside a multigraded Lie algebra (which we denoted $M(E)$, $A(V)$
respectively). Then all the "different" results we obtained are consequences 
of only "one" fact, namely the multigraded Jacobi identity. In the 
line of \cit!{11} it seems to
us that this procedure should be somehow extended to other structures 
defined on a (multigraded) vector space, for example coalgebras, comodules
and then of course to bialgebras such as Hopf algebras and Lie bialgebras.
The latter one was discussed in \cit!{7}.

\heading \totoc \nmb0{5}. Structures and their formal deformations \endheading
 
\subheading{\nmb.{5.1}. Structures}
We fix once for all a $\Bbb Z^{n}$ graded Lie algebra $({\Cal E}, [, ])$ and
a multidegree $\theta \in \Bbb Z^{n}$ such that $\langle \theta,\theta\rangle  + 1 \equiv
0$ (mod 2) (i.e. $\theta$ has one odd number of odd components).
 
By definition, a {\it structure (of degree $\theta$) of\/} ${\Cal E}$ is
an element $P \in {\Cal E}^{\theta}$ such that $[P,P]=0$. We denote by
${\Cal S}^{\theta}({\Cal E})$ the set of structures of degree $\theta$ of
${\Cal E}$.
 
If $P \in {\Cal S}^{\theta}({\Cal E})$, then the adjoint action
$\partial_{P} := \operatorname{ad}P$ of $P$ on ${\Cal E}$ is a differential
homogeneous of degree $\theta$, since 
$\langle \theta,\theta\rangle +1 \equiv 0$ (mod 2). We denote by
$$H({\Cal E}, \partial_{P}) 
     = \bigoplus_{x \in \Bbb Z^{n}} H^{x}({\Cal E},\partial_{P})$$
its cohomology, where
$$H^{x}({\Cal E},\partial_{P}) = {\Cal E}^{x} \cap \ker
     \partial_{P}/\partial_{P} {\Cal E}^{x-\theta}.$$
As $\partial_{P}$ is a derivation of ${\Cal E}$, $H({\Cal
E},\partial_{P})$ has a unique $\Bbb Z^{n}$-graded Lie algebra structure
making the natural map $\ker \partial_{P} \rightarrow H({\Cal E},
\partial_{P})$ a surjective homomorphism of graded Lie algebras.
 
Observe that $[\quad,\quad]$ is a structure of degree $e_{1}=(1,0,\ldots,0)$ of $A({\Cal
E})$ :
$$[\quad,\quad] \in {\Cal S}^{e_{1}} (A({\Cal E})).$$
To avoid confusion as well as to make the notations lighter, we denote
in the sequel by $\Bbb H({\Cal E})$ the space $H(A({\Cal E}),
\partial_{[\quad,\quad]})$ and by $\Bbb D$ the differential $\partial_{[\quad,\quad]}$.
 
As mentioned above, many useful algebraic structures on a vector space
are particular instances of the abstract notion of structure introduced
here (associative algebras, Lie algebras, graded or not, Lie bialgebras
for instance). This leads to a unified way to study these various
algebraic structures, what we shall now illustrate for their formal
deformations.
 
\subheading{\nmb.{5.2}. Formal deformations, Equivalences}
We denote by ${\Cal E}_{\lambda}^{x}$ the space of formal power series in
the parameter $\lambda$ with coefficients in ${\Cal E}^{x}\;(x \in
\Bbb Z^{n})$. The space ${\Cal E}_{\lambda} = \oplus_{x \in \Bbb Z^{n}} {\Cal
E}_{\lambda}^{x}$ has a canonical multigraded Lie algebra structure
extending that of ${\Cal E}$ :
$$\left[\sum_{k} \lambda^{k} X_{k}, \sum_{l} \lambda^{l} Y_{l}\right] =
\sum_{k} \lambda^{k} \sum_{i+j=k} [X_{i}, Y_{j}].$$
 
By definition, a {\it formal deformation} of a structure $P \in {\Cal
S}^{\theta}({\Cal E})$ is an element $P_{\lambda} \in {\Cal S}^{\theta}({\Cal
E}_{\lambda})$ such that $P_{0}=P$. Two such deformations $P_{\lambda}$ and
$P_{\lambda}^{\prime}$ are said to be {\it equivalent\/} if
$P_{\lambda}^{\prime} = \varphi_{\lambda}(P_{\lambda})$ for some
automorphism $\varphi_{\lambda} = \sum_{k} \lambda^{k} \varphi_{k}$
$(\varphi_{k} \in A^{0}({\Cal E}),\; k \in \Bbb N)$ of ${\Cal E}_{\lambda}$
such that $\varphi_{0}= id_{{\Cal E}}$, the identity on ${\Cal E}$. In the
sequel, such a $\varphi_{\lambda}$ will be called an {\it equivalence.}
 
\proclaim{Lemma}
(i)  A mapping $\varphi_{\lambda} :
{\Cal E}_{\lambda} \rightarrow {\Cal E}_{\lambda}$ is
an equivalence if and only if it is solution of a
formal differential equation
$$\frac{d}{d\lambda} \; \varphi_{\lambda} =
\varphi_{\lambda} \circ T_{\lambda},\qquad \varphi_{0} = id_{{\Cal
E}},$$
where $T_{\lambda} \in A^{0}({\Cal E})_{\lambda}$ is a
$\Bbb D$-cocycle.
 
\noindent (ii) Let $\varphi_{\lambda}$ be an
equivalence and let $C_{\lambda} \in A({\Cal
E})_{\lambda}$ be a $\Bbb D$-cocycle. If $\Bbb H^{0}({\Cal
E})=0$, then $\varphi_{\lambda}^{*} C_{\lambda}$ is a
$\Bbb D$-cocycle cohomologous to $C_{\lambda}$.
\endproclaim 
 
Here, $\varphi_{\lambda}^{*}$ denotes the natural
action of $\varphi_{\lambda}$ on $A({\Cal
E})_{\lambda}$ :
$$(\varphi_{\lambda}^{*}C_{\lambda})(X_{0}, \ldots, X_{k}) =
\varphi_{\lambda}(C_{\lambda}(\varphi_{\lambda}^{-1}(X_{0},
\ldots, \varphi_{\lambda}^{-1}(X_{k})).$$
 
\demo{Proof}
(i) Applying $\varphi_{\lambda}^{-1}
\frac{d}{d\lambda}$ to the members of the equation
$$\varphi_{\lambda}([X,Y]) = [\varphi_{\lambda}(X),
\varphi_{\lambda}(Y)]\;\;\;\;(X,Y \in {\Cal E})$$
shows that $T_{\lambda} = \varphi_{\lambda}^{-1}
\frac{d\varphi_{\lambda}}{d\lambda}$ is a
$\Bbb D$-cocycle. Conversely, the unique solution of
$$\frac{d}{d\lambda} \; \varphi_{\lambda} =
\varphi_{\lambda} \circ T_{\lambda}, \;\varphi_{0} =
id_{{\Cal E}},$$
which is given stepwise by
$$(k+1)\varphi_{k+1} = \sum_{i+j=k} \varphi_{i} \circ
T_{j},\; \varphi_{0}=id_{{\Cal E}},$$
is an equivalence if $\Bbb D T_{\lambda} = 0$. Indeed,
as $\varphi_{0} = id_{{\Cal E}}$, it is a bijective
mapping. Moreover,
$$\frac{d}{d\lambda}\;(\varphi_{\lambda}([\varphi_{\lambda}^{-1}
(X), \varphi_{\lambda}^{-1}(Y)])) = \varphi_{\lambda}((\Bbb D
T_{\lambda})(\varphi_{\lambda}^{-1}(X), 
\varphi_{\lambda}^{-1}(Y))) = 0$$
for all $X,Y \in {\Cal E}$. Thus
$$\varphi_{\lambda}([\varphi_{\lambda}^{-1}(X),
\varphi_{\lambda}^{-1}(Y)]) =
\varphi_{\lambda}([\varphi_{\lambda}^{-1}(X),
\varphi_{\lambda}^{-1}(Y)])|_{\lambda =0} = [X,Y].$$ 
(ii)
Assume that $\Bbb D C=0$, where $C \in A({\Cal E})$. As easily
seen, one has 
$$\frac{d}{d\lambda}\;(\varphi_{\lambda}^{*}C) =
\varphi_{\lambda}^{*} [C, \varphi_{\lambda}^{-1}
\frac{d}{d\lambda}\; \varphi_{\lambda}]^{\wedge}.$$ Since
$\Bbb H^{0}({\Cal E}) = 0,\; \varphi_{\lambda}^{-1}
\frac{d}{d\lambda}\;\varphi_{\lambda}$ is a
coboundary. It thus reads $ad\;T_{\lambda}$ for some
$T_{\lambda} \in {\Cal E}_{\lambda}^{0} =
A^{(-1,0)}({\Cal E})$. Noticing that $ad\;T_{\lambda}
= i(T_{\lambda})[\quad,\quad]$, it follows immediately from the
Jacobi identity in $A({\Cal E})$ that
$$[ad\; T_{\lambda}, C]^{\wedge} = i(T_{\lambda}) \Bbb D
C + \Bbb D(i(T_{\lambda})C).$$
Thus
$$\frac{d}{d\lambda}\;(\varphi_{\lambda}^{*}C) =
\varphi_{\lambda}^{*}(\Bbb D i(T_{\lambda})C) =
\Bbb D(\varphi_{\lambda}^{*}(i(T_{\lambda})C))$$
because $\Bbb D C=0$ and, obviously,
$\varphi_{\lambda}^{*} \circ \Bbb D = \Bbb D \circ
\varphi_{\lambda}^{*}$. Therefore,
$$\varphi_{\lambda}^{*}C = C + \Bbb D \int_{0}^{\lambda}
\varphi_{\mu}^{*}(i(T_{\mu})C)d\mu.$$
Now, if each component $C_{k}$ of $C_{\lambda}$ is a
$\Bbb D$-cocycle, then
$$\varphi_{\lambda}^{*}C_{\lambda} = C_{\lambda} +
\Bbb D(\sum_{k} \lambda^{k} \int_{0}^{\lambda}
\varphi_{\mu}^{*}(i(T_{\mu})C_{k})d\mu)$$
is cohomologous to $C_{\lambda}$. 
\qed\enddemo 
 
\proclaim{\nmb.{5.3}. Proposition}
Let ${\Cal E}$ be a Lie subalgebra of $A(V)$ for some
$(n-1)$-graded vector space $V$. Assume that ${\Cal
E}^{(-1,*)} = V^{*}$. Then
 
(i) If $C \in A^{0}({\Cal E})$ is a $\Bbb D$-cocycle,
then $C=\Bbb D T$ for some $T \in A^{0}(V)$ such that
$[T,{\Cal E}] \subset {\Cal E}$.
 
(ii) The equivalences $\varphi_{\lambda}$ of ${\Cal E}$
are the mappings of the form $S_{\lambda}^{*}$ where
$S_{\lambda} \in A^{0}(V)_{\lambda}$ and $S_{0}=
id_{V}$.
\endproclaim
 
\demo{Proof} (i) Set $T = C|V$. Then $T \in A^{0}(V)$ and
$C=\Bbb D T$. Indeed, let $Y \in A^{(k,y)}(V)$. For
$k=-1$, $C(Y) = (\Bbb D T)(Y)$ by definition of $T$. Now,
by induction on $k$, if $X \in A^{(-1,x)}(V) = V^{x}$,
then
$$\align
(-1)^{\langle x,y\rangle } i(X)C(Y) &= (\Bbb D C)(X,Y) +
C([X,Y])+(-1)^{k+\langle x,y\rangle }[Y,C(X)] \\
&= [[X,Y],T] + (-1)^{k+\langle x,y\rangle } [Y,[X,T]] \\
&= [X,(\Bbb D T)(Y)] = (-1)^{\langle x,y\rangle }i(X)((\Bbb D T)(Y)).
\endalign$$
Thus $C(Y) = (\Bbb D T)(Y)$ for all $Y$.
 
(ii) It is clear that $S_{\lambda}^{*}$ is an
equivalence. Conversely, if $\varphi_{\lambda}$ is an
equivalence, we know that $\varphi_{\lambda}^{-1}
\frac{d}{d\lambda}\; \varphi_{\lambda}$ is a
$\Bbb D$-cocycle. It is thus of the form $\Bbb D
T_{\lambda}$ for some $T_{\lambda} \in
A^{0}(V)_{\lambda}$. The equation
$$\frac{d}{d\lambda} \; S_{\lambda} = S_{\lambda}
\circ T_{\lambda},\; S_{0} = id_{V},$$
has, obviously, a unique solution. But then, for an
arbitrary $X \in {\Cal E}$, one has
$$\frac{d}{d\lambda}\; \varphi_{\lambda}(X) =
\varphi_{\lambda}(\Bbb D T_{\lambda})(X),
\frac{d}{d\lambda}\; S_{\lambda}^{*} X =
S_{\lambda}^{*} (\Bbb D T_{\lambda})(X)$$
and thus $\varphi_{\lambda}$ and $S_{\lambda}^{*}$
coincide on ${\Cal E}$ since $\varphi_{0} =
S_{0}^{*}$. \qed\enddemo

\subheading{\nmb.{5.4}} We now turn to generalize to
arbitrary structures $P \in {\Cal S}^{\theta}({\Cal
E})$ the results obtained for the Lie algebras in
\cit!{1}, \cit!{6}. We only indicate the non
obvious adaptations of the proofs, referring otherwise
the reader to the appropriate papers. As before,
$({\Cal E}, [\quad,\quad])$ denotes a $\Bbb Z^{n}$-graded Lie
algebra and $\theta \in \Bbb Z^{n}$ is assumed to be
such that $\langle \theta,\theta\rangle  +1 \equiv 0$ (mod 2). We
also denote by $Pol({\Cal E})$ the space of polynomials
on ${\Cal E}$. Let $\eta$ be the map $A({\Cal E})
\rightarrow Pol({\Cal E})$ given by
$$\eta(C) : X \rightarrow \frac{1}{(k+1)!} \; C(X,
\ldots, X)$$
for $C \in A^{(k,c)}({\Cal E})$ and set
$$\eta_{P}(C) = \eta(C)(P)$$
for $P \in {\Cal S}^{\theta}({\Cal E})$.
 
\proclaim{\nmb.{5.5}. Lemma}
Let $C \in A^{(k,c)}({\Cal E})$, $P \in {\Cal
S}^{\theta}({\Cal E})$ and $X \in {\Cal E}^{x}$ be given
such that $k >0$ and $\langle x,x\rangle +1 \equiv 0$ (mod 2).
Then we have
$$\gather 
\eta(\Bbb D C)(X) = (-1)^{\langle x,c\rangle }[X,\eta(C)(X)] -
     \frac{1}{2}\;\eta(i[X,X])C)(X) \\
\partial_{P} \eta_{P}(i(X)C)+\eta_{P}(i(\partial_{P}X)C)
     =(-1)^{\langle \theta,x+c\rangle }[X,\eta_{P}(C)]
     \text{ if } \Bbb D C=0.
\endgather$$
\endproclaim
 
\demo{Proof} This we get by straightforward computations
applying theorem \nmb!{4.4} to expand $\Bbb D C(X,\ldots,X)$ and
$\Bbb D C(X,P,\ldots, P)$ respectively. \qed\enddemo
 
Now comes our main result about formal deformation of
structures of ${\Cal E}$.
 
\proclaim{\nmb.{5.6}. Theorem}
(i) Let $C_{\lambda} \in (\bigoplus_{k\theta+c=0}
A^{(k,c)}({\Cal E}))_{\lambda}$ be such that
$\Bbb D C_{\lambda}=0$. For each $P \in {\Cal
S}^{\theta}({\Cal E})$, the unique solution of the
equation
$$\frac{d}{d\lambda}\;P_{\lambda} +
\eta(C_{\lambda})(P_{\lambda})=0,\;P_{0}=P,$$
is a formal deformation of $P$, which will be said to be
associated to the cocycle $C_{\lambda}$.
 
(ii) Formal deformations which are associated to cohomologous
cocycles are equivalent.
 
(iii) If $\Bbb H^{0}({\Cal E})=0$ and $P_{\lambda}$ is
associated to $C_{\lambda}$, then each deformation
equivalent to $P_{\lambda}$ is associated to a cocycle
cohomologous to $C_{\lambda}$.
 
(iv) For a given $P \in {\Cal S}^{0}({\Cal E})$, the
image of
$$\eta_{P\#} : \Bbb H({\Cal E}) \rightarrow {\Cal E}$$
lies in the center of $H(E, \partial_{P})$. If
$H^{\theta}({\Cal E},\partial_{P}) \subset im\;\eta_{P\#}$,
then each deformation of $P$ is associated to some cocycle.
\endproclaim
 
\demo{Proof} (i) The proof goes as in (\cit!{1},
Proposition 15.2), without major change : simply, substitute
the first equation of lemma \nmb!{5.5} to Proposition 15.1 in
\cit!{1}.
 
(ii) Assume that $C_{\lambda}^{\prime} = C_{\lambda} +
\Bbb D A_{\lambda}$ and denote by $P_{\lambda}$ and
$P_{\lambda}^{\prime}$ the deformations of the same $P
\in {\Cal S}^{\theta}({\Cal E})$ associated to
$C_{\lambda}$ and $C_{\lambda}^{\prime}$ respectively.
Set 
$$B_{\lambda} = \sum_{k} \lambda^{k}
\int_{0}^{\lambda}
(\varphi_{\mu}^{*}(i(T_{\mu})C_{k})d\mu$$
(see the end of the proof of the Lemma in \nmb!{5.2}). Then,
the equations
$$\frac{d}{d\lambda}\;\varphi_{\lambda} = 
\varphi_{\lambda} \circ ad\; T_{\lambda},\;
\varphi_{0} = id_{{\Cal E}},$$
and
$$T_{\lambda} =
\varphi_{\lambda}^{-1}(\eta(A_{\lambda}-B_{\lambda})
(P_{\lambda}^{\prime}))$$ have a unique solution
$\varphi_{\lambda}, T_{\lambda}$, where $\varphi_{\lambda}$ is
an equivalence and $T_{\lambda} \in A^{(-1,0)}({\Cal
E})_{\lambda}$. Indeed, it follows from the first that
$\varphi_{k}$ is uniquely expressed in terms of
$T_{0}, \ldots, T_{k-1}$. The same holds true for the
$k$-th component of $\varphi_{\lambda}^{-1}$, a
polynomial in $\varphi_{0}, \ldots, \varphi_{k}$.
Thus, the $k$-th component of the right member of the
second equation only depends on $T_{0}, \ldots,
T_{k-1}$. It follows that the two equations may be
uniquely solved by induction. Taking account of the
fact that
$$\varphi_{\lambda}^{*} C_{\lambda} =
C_{\lambda}^{\prime} + \Bbb D(B_{\lambda}-A_{\lambda}),$$
one easily sees that $P_{\lambda}^{\prime}$ and
$\varphi_{\lambda}(P_{\lambda})$ both are associated
to $C_{\lambda}^{\prime}$. As $P_{0}^{\prime} =
\varphi_{0}(P_{0})=P$, one has thus
$P_{\lambda}^{\prime} =
\varphi_{\lambda}(P_{\lambda})$.
 
(iii) Let $\varphi_{\lambda}$ be an equivalence. As
$\Bbb H^{0}({\Cal E}) = 0$, one has
$\frac{d\varphi_{\lambda}}{d\lambda} =
\varphi_{\lambda} \circ ad\;T_{\lambda}$ for some
$T_{\lambda} \in {\Cal E}_{\lambda}^{0}$, and by
lemma in \nmb!{5.2},
$\varphi_{\lambda}^{*}C_{\lambda}=C_{\lambda} + \Bbb D
B_{\lambda}$. Computing
$\frac{d}{d\lambda}\;\varphi_{\lambda} (P_{\lambda})$
easily shows that $\varphi_{\lambda}(P_{\lambda})$ is
associated to
$C_{\lambda}+\Bbb D(B_{\lambda}+\varphi_{\lambda}(T_{\lambda}))$.
 
(iv) The fact that $im\;\eta_{P\#}$ lies in the center
of $H({\Cal E}, \partial_{P})$ follows immediately
from the second equation in lemma \nmb!{5.5}. The proof of the
second part of (iv) goes as in (\cit!{6}, Prop.4.4).
\qed\enddemo

\enddocument